\documentclass[11pt]{article}

\usepackage{amssymb}
\usepackage{amsfonts}
\usepackage{amsmath}
\usepackage{amsthm}
\usepackage[margin=1.00in]{geometry}
\usepackage{enumerate}
\usepackage{amstext}
\usepackage{layout}

\numberwithin{equation}{section}

\newtheorem{theorem}{Theorem}[section]
\newtheorem{corollary}{Corollary}[section]
\newtheorem{lemma}{Lemma}[section]

\newtheorem{remark}{Remark}[section]
\newtheorem{example}{Example}[section]

\begin{document}

\noindent {\bf\large{A refinement of the
Kolmogorov-Marcinkiewicz-Zygmund strong law of large numbers}}

\vskip 0.3cm

\noindent {\bf Deli Li\footnote{Deli Li, Department of Mathematical
Sciences, Lakehead University, Thunder Bay, Ontario, Canada P7B
5E1\\ e-mail: dli@lakeheadu.ca} $\cdot$ Yongcheng
Qi\footnote{Department of Mathematics and
Statistics, University of Minnesota Duluth, Duluth, Minnesota 55812, USA\\
e-mail: yqi@d.umn.edu}$\cdot$ Andrew Rosalsky\footnote{Andrew
Rosalsky, Department of Statistics, University of Florida,
Gainesville, Florida 32611, U.S.A.\\ e-mail: rosalsky@stat.ufl.edu}
\footnote{Corresponding author: Andrew Rosalsky (Telephone:
1-352-273-2983, FAX: 1-352-392-5175)}}

\vskip 0.3cm

\noindent {\bf Abstract} Let $ \{X_{n};~n \geq 1 \}$ be a sequence
of independent copies of a real-valued random variable $X$ and set
$S_{n} = X_{1} + \cdots + X_{n}, ~n \geq 1$. This paper is devoted
to a refinement of the classical Kolmogorov-Marcinkiewicz-Zygmund
strong law of large numbers. We show that for $0 < p < 2$,
\[
\sum_{n = 1}^{\infty} \frac{1}{n}
\left(\frac{|S_{n}|}{n^{1/p}}\right) < \infty~~\mbox{almost surely}
\]
if and only if
\[
\left \{
\begin{array}{ll}
\mbox{$\displaystyle \mathbb{E}|X|^{p} < \infty, $}
& \mbox{if~~$\displaystyle 0 < p < 1$}\\
&\\
\mbox{$\displaystyle \mathbb{E}X = 0,~\sum_{n=1}^{\infty}
\frac{\left|\mathbb{E}XI\{|X| \leq n\} \right|}{n} <
\infty,~\mbox{and}~\sum_{n=1}^{\infty} \frac{\int_{\min\{u_{n},
n\}}^{n} \mathbb{P}(|X| > t) dt}{n} < \infty,$} & \mbox{if
~~$\displaystyle p = 1$}\\
&\\
\mbox{$\displaystyle \mathbb{E}X = 0 ~\mbox{and}~\int_{0}^{\infty}
\mathbb{P}^{1/p}(|X|
> t) dt < \infty,$} & \mbox{if~~$\displaystyle 1 < p < 2,$}
\end{array}
\right.
\]
where $u_{n} = \inf \left\{t:~ \mathbb{P}(|X| > t) < \frac{1}{n}
\right\}, ~n \geq 1$. Versions of above results in a Banach space
setting are also presented. To establish these results, we invoke
the remarkable Hoffmann-J{\o}rgensen (1974) inequality to obtain
some general results for sums of the form
$\sum_{n=1}^{\infty}a_{n} \left\|\sum_{i=1}^{n}V_{i}\right\|$
(where $\{V_{n}; ~n \geq 1\}$ is a sequence of independent Banach space
valued random variables and $a_{n} \geq 0, ~n \geq 1$) which may be of
independent interest but which we apply to $\sum_{n = 1}^{\infty} \frac{1}{n}
\left(\frac{|S_{n}|}{n^{1/p}}\right)$.

~\\

\noindent {\bf Keywords} Kolmogorov-Marcinkiewicz-Zygmund strong law
of large numbers $\cdot$ Sums of i.i.d. random variables $\cdot$
Real separable Banach space $\cdot$ Rademacher type $p$ Banach space
$\cdot$ Stable type $p$ Banach space

\vskip 0.3cm

\noindent {\bf Mathematics Subject Classification (2000)} Primary:
60F15; Secondary: 60B12 $\cdot$ 60G50

\vskip 0.3cm

\noindent {\bf Running Head}: Strong law of large numbers

\section{Introduction and preliminaries}

Throughout, let $(\mathbf{B}, \| \cdot \| )$ be a real separable
Banach space equipped with its Borel $\sigma$-algebra $\mathcal{B}$
($=$ the $\sigma$-algebra generated by the class of open subsets of
$\mathbf{B}$ determined by $\|\cdot\|$) and let $ \{X_{n};~n \geq 1
\}$ be a sequence of independent copies of a {\bf B}-valued random
variable $X$ defined on a probability space $(\Omega, \mathcal{F},
\mathbb{P})$. As usual, let $S_{n} = \sum_{k=1}^{n} X_{k},~ n \geq
1$ denote their partial sums. If $0 < p < 2$ and if $X$ is a
real-valued random variable (that is, if $\mathbf{B} = \mathbb{R}$),
then
\[
\lim_{n \rightarrow \infty} \frac{S_{n}}{n^{1/p}} = 0 ~~\mbox{almost
surely (a.s.)}
\]
if and only if
\[
\mathbb{E}|X|^{p} < \infty ~~\mbox{where}~~ \mathbb{E}X = 0
~~\mbox{whenever}~~ p \geq 1.
\]
This is the celebrated Kolmogorov-Marcinkiewicz-Zygmund strong law
of large numbers (SLLN); see Kolmogoroff [8] for $p = 1$ and
Marcinkiewicz and Zygmund [11] for $p \neq 1$.

The classical Kolmogorov SLLN in real separable Banach spaces was
established by Mourier [14]. The extension of the
Kolmogorov-Marcinkiewicz-Zygmund SLLN to $\mathbf{B}$-valued random
variables is independently due to Azlarov and Volodin [1] and de
Acosta [3].

\vskip 0.3cm

\begin{theorem}
{\rm (Azlarov and Volodin [1] and de Acosta [3])}. Let $0 < p < 2$
and let $\{X_{n}; ~n \geq 1\}$ be a sequence of independent copies
of a $\mathbf{B}$-valued random variable $X$. Then
\[
\lim_{n \rightarrow \infty} \frac{S_{n}}{n^{1/p}} = 0~~\mbox{a.s.}
\]
if and only if
\[
\mathbb{E}\|X\|^{p} < \infty~~\mbox{and}~~\frac{S_{n}}{n^{1/p}}
\rightarrow_{\mathbb{P}} 0.
\]
\end{theorem}

\vskip 0.3cm

De Acosta [3] also provides a remarkable characterization of
Rademacher type $p$ Banach spaces. (Technical definitions such as
$\mathbf{B}$ being of Rademacher type $p$ will be reviewed below.)
Specifically, de Acosta [3] proved the following theorem.

\vskip 0.3cm

\begin{theorem}
{\rm (de Acosta [3])}. Let $1 \leq p < 2$. Then the following two
statements are equivalent:
\begin{align*}
& {\bf (i)} \quad \mbox{The Banach space $\mathbf{B}$ is of
Rademacher type $p$.}\\
& {\bf (ii)} \quad \mbox{For every sequence $\{X_{n}; ~n \geq 1 \}$
of independent copies of a $\mathbf{B}$-valued variable $X$},
\end{align*}
\[
\lim_{n \rightarrow \infty} \frac{S_{n}}{n^{1/p}} = 0~~\mbox{a.s. if
and only if}~~\mathbb{E}\|X\|^{p} < \infty~~\mbox{and}~~\mathbb{E}X
= 0.
\]
\end{theorem}

\vskip 0.3cm

At the origin of the current investigation is the following recent
and striking result by Hechner and Heinkel [4].

\vskip 0.3cm

\begin{theorem}
{\rm (Hechner and Heinkel [4])}. Suppose that $\mathbf{B}$ is of
stable type $p$ ($1 < p < 2$) and let $\{X_{n}; ~n \geq 1 \}$ be a
sequence of independent copies of a $\mathbf{B}$-valued variable $X$
with $\mathbb{E}X = 0$. Then
\begin{equation}
\sum_{n=1}^{\infty} \frac{1}{n}
\left(\frac{\mathbb{E}\|S_{n}\|}{n^{1/p}} \right) < \infty
\end{equation}
if and only if
\begin{equation}
\int_{0}^{\infty} \mathbb{P}^{1/p}(\|X\| > t) dt < \infty.
\end{equation}
\end{theorem}

\vskip 0.3cm

This result of Hechner and Heinkel [4] is new even in the case where
the Banach space $\mathbf{B}$ is the real line. We note that,
under the hypotheses of Theorem 1.3, (1.1) and (1.2) each imply that
\[
\lim_{n \rightarrow \infty} \frac{\mathbb{E}\|S_{n}\|}{n^{1/p}} = 0.
\]
This follows from Theorem 1.3, Remark 2.2 below, Theorem 1.2, the fact
that $\mathbf{B}$ is of Rademacher type $p$ (see the discussion below),
and Theorem 2 of Korzeniowski [9]. This result of Korzeniowski extends
to a Banach space setting Corollary 12 of Klass [7].

Inspired by the above discovery by Hechner and Heinkel [4], in the
current work we obtain sets of necessary and sufficient conditions
for
\begin{equation}
\sum_{n=1}^{\infty} \frac{1}{n} \left(\frac{\|S_{n}\|}{n^{1/p}}
\right) < \infty~~\mbox{a.s.}
\end{equation}
for the three cases: $0 < p < 1$, $p = 1$, $1 < p < 2$. Moreover, we
obtain necessary and sufficient conditions for
\begin{equation}
\sum_{n=1}^{\infty} \frac{1}{n} \left(\frac{\mathbb{E}\|S_{n}\|}{n}
\right) < \infty.
\end{equation}
Again, these results are new when $\mathbf{B} = \mathbb{R}$; see
Theorem 2.5. The current work complements the investigation by
Hechner and Heinkel [4].

While it is immediate that (1.1) implies (1.3), we will see that
(1.1) and (1.3) are indeed equivalent if $1 < p < 2$ (Theorem 2.1).
However, there is a gap between the cases $0< p \leq 1$ and $1 < p <
2$ as we will see that if $0 < p \leq 1$, then (1.4) and (1.3) are
not equivalent (Theorems 2.2 and 2.4).  Moreover, Theorems 2.2 and
2.4 also reveal that there is a gap between the cases $p = 1$ and $0
< p < 1$.

The most delicate case is that for $p = 1$. We show assuming
$\mathbf{B}$ is of stable type $1$ that (1.3) holds if and only if
the conditions (2.11), (2.12), and (2.13) are fulfilled (Theorem
2.3).

We now review various technical definitions pertaining to a
$\mathbf{B}$-valued random variable $X$ or to the Banach space
$\mathbf{B}$ itself.

The {\it expected value} or {\it mean} of $X$, denoted
$\mathbb{E}X$, is defined to be the {\it Pettis integral} provided
it exists. That is, $X$ has expected value $\mathbb{E}X \in
\mathbf{B}$ if $\varphi(\mathbb{E}X) = \mathbb{E}\varphi(X)$ for
every $\varphi \in \mathbf{B}^{*}$ where $\mathbf{B}^{*}$ denotes
the {\it (dual)} space of all continuous linear functionals on
$\mathbf{B}$. If $\mathbb{E}\|X\| < \infty$, then (see, e.g., Taylor
[17, p. 40]) $X$ has an expected value. But the expected value can
exist when $\mathbb{E}\|X\| = \infty$. For an example, see Taylor
[17, p. 41].

Let $\{R_{n}; ~n \geq 1\}$ be a {\it Rademacher sequence}; that is,
$\{R_{n};~n \geq 1\}$ is a sequence of independent and identically
distributed (i.i.d.) random variables with $\mathbb{P}\left(R_{1} =
1\right) = \mathbb{P}\left(R_{1} = -1\right) = 1/2$. Let
$\mathbf{B}^{\infty} =
\mathbf{B}\times\mathbf{B}\times\mathbf{B}\times \cdots$ and define
\[
\mathcal{C}(\mathbf{B}) = \left\{(v_{1}, v_{2}, ...) \in
\mathbf{B}^{\infty}: ~\sum_{n=1}^{\infty} R_{n}v_{n}
~~\mbox{converges in probability} \right\}.
\]
Let $1 \leq p \leq 2$. Then $\mathbf{B}$ is said to be of {\it
Rademacher type $p$} if there exists a constant $0 < C < \infty$
such that
\[
\mathbb{E}\left\|\sum_{n=1}^{\infty}R_{n}v_{n} \right\|^{p} \leq C
\sum_{n=1}^{\infty}\|v_{n}\|^{p}~~\mbox{for all}~(v_{1}, v_{2}, ...)
\in \mathcal{C}(\mathbf{B}).
\]
Hoffmann-J{\o}rgensen and Pisier [6] proved for $1 \leq p \leq 2$
that $\mathbf{B}$ is of Rademacher type $p$ if and only if there
exists a constant $0 < C < \infty$ such that
\[
\mathbb{E}\left\|\sum_{k=1}^{n}V_{k} \right\|^{p} \leq C
\sum_{k=1}^{n} \mathbb{E}\left\|V_{k}\right\|^{p}
\]
for every finite collection $\{V_{1}, ..., V_{n} \}$ of independent
mean $0$ $\mathbf{B}$-valued random variables.

If $\mathbf{B}$ is of Rademacher type $p$ for some $p \in (1, 2]$,
then it is of Rademacher type $q$ for all $q \in [1, p)$. Every real
separable Banach spaces is of Rademacher type (at least) $1$.

Let $0 < p \leq 2$ and let $\{\Theta_{n}; ~n \geq 1 \}$ be a
sequence of i.i.d. stable random variables each with characteristic
function $\psi(t) = \exp \left\{-|t|^{p}\right\}, ~- \infty < t <
\infty$. Then $\mathbf{B}$ is said to be of {\it stable type $p$} if
$\sum_{n=1}^{\infty} \Theta_{n}v_{n}$ converges a.s. whenever
$\{v_{n}: ~n \geq 1\} \subseteq \mathbf{B}$ with
$\sum_{n=1}^{\infty} \|v_{n}\|^{p} < \infty$. Equivalent
characterizations of a Banach space being of stable type $p$,
properties of stable type $p$ Banach spaces, as well as various
relationships between the conditions ``Rademacher type $p$" and
``stable type $p$" may be found in Maurey and Pisier [12],
Woyczy\'{n}ski [18], Marcus and Woyczy\'{n}ski [13], Rosi\'{n}ski
[16], Pisier [15], and Ledoux and Talagrand [10]. Some of these
properties and relationships will now be summarized:

\begin{description}
\item
{\bf (1)} Every real separable Banach space $\mathbf{B}$ is of
stable type $p$ for all $p \in (0, 1)$.

\item
{\bf (2)} For $1 \leq p < 2$, $\mathbf{B}$ is of stable type $p$ if
and only if $\mathbf{B}$ is of Rademacher type $p_{1}$ for some
$p_{1} \in (p, 2]$.

\item
{\bf (3)} $\mathbf{B}$ is of stable type $2$ if and only if
$\mathbf{B}$ is of Rademacher type $2$.
\end{description}

\noindent Consequently:

\begin{description}
\item
{\bf (4)} If $\mathbf{B}$ is of stable type $p$ for some $p \in [1,
2]$, then $\mathbf{B}$ is of Rademacher type $p$.

\item
{\bf (5)} If $\mathbf{B}$ is of stable type $p$ for some $p \in [1,
2]$, then $\mathbf{B}$ is of stable type $q$ for all $q \in (0, p)$.

\item
{\bf (6)} If $\mathbf{B}$ is of stable type $p$ for some $p \in [1,
2)$, then $\mathbf{B}$ is of stable type $q$ for some $q \in (p,
2]$.
\end{description}
The property {\bf (6)} is the Maurey-Pisier [12] theorem.

For $q \geq 2$, the $\mathcal{L}_{q}$-spaces and $\ell_{q}$-spaces
are of stable type $2$ while for $1 \leq q < 2$, the
$\mathcal{L}_{q}$-spaces and $\ell_{q}$-spaces are of Rademacher
type $q$ and are of stable type $p$ for all $p \in (0, q)$ but are
not of stable type $q$. Every real separable Hilbert space and real
separable finite dimensional Banach space is of stable type $2$. In
particular, the real line is of stable type $2$.

The plan of the paper is as follows. The main results are stated in
Section 2. Some general results for
$\sum_{n=1}^{\infty}a_{n}\|\sum_{k=1}^{n}V_{k} \|$ (where the $a_{n}
\geq 0$ and $\{V_{k}; ~k \geq 1 \}$ is a sequence of independent
$\mathbf{B}$-valued random variables) is established in Section 3;
these results are key components in the proofs of the main results.
The main results are proved in Sections 4 and 5.

Finally, the symbol $C$ denotes throughout a generic constant ($0 <
C < \infty$) which is not necessarily the same one in each
appearance.

\section{Statement of the main results}

With the preliminaries accounted for, the main results may be
stated. We begin with the case where $1 < p < 2$.

\vskip 0.3cm

\begin{theorem}
Let $ \{X_{n};~n \geq 1 \}$ be a sequence of independent copies of a
$\mathbf{B}$-valued random variable $X$. Let $1 < p < 2$. Then
\begin{equation}
\sum_{n = 1}^{\infty} \frac{1}{n} \left(
\frac{\mathbb{E}\|S_{n}\|}{n^{1/p}}\right) < \infty
\end{equation}
if and only if
\begin{equation}
\sum_{n = 1}^{\infty} \frac{1}{n} \left(
\frac{\|S_{n}\|}{n^{1/p}}\right) < \infty~~\mbox{a.s.}
\end{equation}
Furthermore, each of (2.1) and (2.2) implies that
\begin{equation}
\mathbb{E}X = 0~~\mbox{and}~~\int_{0}^{\infty}
\mathbb{P}^{1/p}(\|X\|
> t) dt < \infty
\end{equation}
and
\begin{equation}
\lim_{n \rightarrow \infty} \frac{S_{n}}{n^{1/p}} = 0~~{a.s.}
\end{equation}
\end{theorem}

\vskip 0.3cm

\begin{remark}
The conditions (2.3) and (2.4) do not necessarily imply that (2.1)
or (2.2) hold. A counterexample provided by Hechner and Heinkel [4,
Section 5] is such that (2.3) and (2.4) hold but (2.1) fails. Since
(2.1) and (2.2) are equivalent, (2.2) also fails for that
counterexample. The counterexample pertained to the Banach space
$\ell_{p}$ $(1 < p < 2)$.
\end{remark}

\vskip 0.3cm

\begin{remark}
For $1 < p < 2$, the second half of (2.3) implies that
$\mathbb{E}\|X\|^{p} < \infty$ as was noted by Hechner and Heinkel
[4]. But the converse implication is false. To see this, suppose
that
\[
\mathbb{P}(\|X\| > t) = \frac{e^{ep + 1}}{t^{p} (\ln t) (\ln \ln
t)^{\alpha}}, ~ t \geq e^{e}
\]
where $\alpha > 1$. Then $\mathbb{E} \|X\|^{p} < \infty$ but
\[
\int_{0}^{\infty} \mathbb{P}^{1/p}(\|X\| > t) dt = C +
\int_{e^{e}}^{\infty} \frac{e^{e + p^{-1}}}{t (\ln t)^{1/p}(\ln \ln
t)^{\alpha/p}}dt = \infty.
\]
\end{remark}

\vskip 0.3cm

Combining Theorem 2.1 above and Theorem 5 of Hechner and Heinkel
[4], we immediately obtain the following result.

\vskip 0.3cm

\begin{corollary}
Let $X$ be a {\bf B}-valued random variable and let $1 < p < 2$. If
$\mathbf{B}$ is of stable type $p$, then (2.1), (2.2), and (2.3) are
equivalent.
\end{corollary}

\vskip 0.3cm

\begin{remark}
The example of Hechner and Heinkel [4, Section 5] referred to in
Remark 2.1 above also shows that in Corollary 2.1, the stable type
$p$ ($1 < p < 2$) hypothesis cannot be weakened to a Rademacher type
$p$ ($1 < p < 2$) hypothesis.
\end{remark}

\vskip 0.3cm

\begin{remark}
An example of Hechner and Heinkel [4, Example 1] shows that
Corollary 2.1 can fail if the condition (2.3) is weakened to
\[
\mathbb{E}X = 0 ~~\mbox{and}~~\mathbb{E}\|X\|^{p} < \infty.
\]
\end{remark}

\vskip 0.3cm

We now consider the case where $p = 1$. In general (2.1) and (2.2)
are not equivalent when $p = 1$; see Remark 2.7 below. But we have
the following result.

\vskip 0.3cm

\begin{theorem}
Let $X$ be a {\bf B}-valued random variable $X$. Then
\begin{equation}
\sum_{n = 1}^{\infty} \frac{1}{n}\left(\frac{\mathbb{E}\|S_{n}\|}{n}
\right) < \infty
\end{equation}
if and only if
\begin{equation}
\sum_{n = 1}^{\infty} \frac{1}{n} \left(\frac{\|S_{n}\|}{n} \right)
< \infty~~\mbox{a.s.}
\end{equation}
and
\begin{equation}
\mathbb{E}X = 0~~\mbox{and}~~\mathbb{E}\|X\| \ln (1 + \|X\|) <
\infty.
\end{equation}
\end{theorem}

\vskip 0.3cm

In the next theorem, we provide necessary and sufficient conditions
for (2.6) to hold assuming that $\mathbf{B}$ is of stable type $1$.
It is the most delicate result in this paper.

Let $X$ be a {\bf B}-valued random variable. For each $n \geq 1$, we
define the quantile $u_{n}$ of order $\left(1 - \frac{1}{n}\right)$
of $\|X\|$ as follows:
\[
u_{n} = \inf\left\{t:~ \mathbb{P}(\|X\| \leq t) > 1 - \frac{1}{n}
\right\} = \inf \left\{t:~ \mathbb{P}(\|X\| > t) < \frac{1}{n}
\right\}.
\]
If
\begin{equation}
\mathbb{E}\|X\| < \infty,
\end{equation}
then it is easy to show that
\begin{equation}
\lim_{n \rightarrow \infty} \frac{u_{n}}{n} = 0.
\end{equation}

\vskip 0.3cm

\begin{theorem}
Let $\mathbf{B}$ be a Banach space of stable type $1$. Let
$\{X_{n};~n \geq 1 \}$ be a sequence of independent copies of a
$\mathbf{B}$-valued random variable $X$. Then
\begin{equation}
\sum_{n=1}^{\infty} \frac{1}{n}\left(\frac{\|S_{n}\|}{n} \right) <
\infty~~\mbox{a.s.}
\end{equation}
if and only if following three conditions are fulfilled:
\begin{equation}
\mathbb{E}\|X\| < \infty~~\mbox{and}~~\mathbb{E}X = 0;
\end{equation}
\begin{equation}
\sum_{n=1}^{\infty} \frac{\left\|\mathbb{E}XI\{\|X\| \leq n\}
\right\|}{n} < \infty;
\end{equation}
\begin{equation}
\sum_{n=1}^{\infty} \frac{\int_{\min\{u_{n}, n\}}^{n}
\mathbb{P}(\|X\| > t) dt}{n} < \infty.
\end{equation}
\end{theorem}

\vskip 0.3cm

\begin{remark}
Set
\[
G(t) = \mathbb{P}(\|X\| > t), ~~\ell(t) = \sum_{t \leq n < 1/G(t)} \frac{1}{n},
~~\mbox{and}~ \pi(t) = \sum_{n=1}^{\infty} \frac{1}{n} I_{(u_{n}, n]}(t),
~~t \geq 0.
\]
Since
\[
u_{n} < t \leq n \Longrightarrow t \leq n < \frac{1}{G(t)}
\Longrightarrow u_{n} \leq t \leq n,
\]
we have $\pi(t) \leq \ell(t)$ for all $t \geq 0$ with equality holding if $t \neq u_{n}$ for all
$n \geq 1$. Hence we have
\[
\int_{0}^{\infty} \pi(t) G(t)dt = \sum_{n=1}^{\infty} \frac{1}{n} \int_{\min\{u_{n}, n\}}^{n} G(t)dt
= \int_{0}^{\infty}\ell(t) G(t)dt
\]
and so we see that (2.13) is equivalent to
\[
\int_{0}^{\infty}\ell(t)G(t) dt < \infty.
\]
An elementary computation shows that
\[
\ell(t) = \ln^{+} \frac{1}{Q(t)} + \gamma(t)
~~\mbox{where}~~Q(t) = (t + 1)G(t)~~\mbox{and}~~|\gamma(t)| \leq \frac{1}{t+1}
~~\mbox{for all}~ t \geq 0.
\]
Hence, if $\mathbb{E}\|X\| < \infty$, we see that (2.13) is equivalent to
\begin{equation}
\int_{0}^{\infty} G(t) \ln^{+} \frac{1}{Q(t)} dt < \infty.
\end{equation}
For instance, if $G(t) \leq C(t \ln t)^{-1}(\ln \ln t)^{-\beta}$ for $t \geq 16$
where $\beta > 2$, then (2.14) and (2.13) hold but we may have
$\mathbb{E}\|X\| \ln^{\delta}(1 + \|X\|) = \infty$ for all $\delta > 0$.
\end{remark}

Combining Theorems 2.2 and 2.3, we obtain the following result.

\vskip 0.3cm

\begin{corollary}
Let $X$ be a $\mathbf{B}$-valued random variable. If $\mathbf{B}$ is
of stable type $1$, then (2.5) and (2.7) are equivalent.
\end{corollary}

\vskip 0.3cm

Another corollary of Theorem 2.3 is as follows.

\vskip 0.3cm

\begin{corollary}
Suppose $\mathbf{B}$ is of stable type $1$ and let $X$ be a
symmetric $\mathbf{B}$-valued random variable. Then (2.6) holds if
\begin{equation}
\mathbb{E}\|X\| \ln^{\delta} (1 + \|X\|) < \infty~~\mbox{for some}~
\delta > 0.
\end{equation}
\end{corollary}

\vskip 0.3cm

\begin{remark}
It follows from Corollary 2.3 that the moment condition (2.7) in the
implication ((2.7) $\Rightarrow$ (2.6)) which is immediate from
Corollary 2.2 can be weakened to (2.15) if it is assumed that $X$ is
symmetric.
\end{remark}

\vskip 0.3cm

\begin{remark}
When $p = 1$, (2.1) and (2.2) are not equivalent. To see this, let
$\{X_{n}; ~ n \geq 1 \}$ be a sequence of independent copies of a
symmetric real-valued random variable $X$ with
\[
\mathbb{E}|X| \ln (1 + |X|) = \infty ~~\mbox{but}~~\mathbb{E}|X|
\ln^{\delta}(1 + |X|) < \infty ~~\mbox{for some}~ 0 < \delta < 1.
\]
Then by Corollary 2.3, (2.2) holds with $p = 1$ but by Corollary 2.2
(or by Theorem 2.2), (2.1) fails with $p = 1$.
\end{remark}

\vskip 0.3cm

We now consider the case where $0 < p < 1$.

\vskip 0.3cm

\begin{theorem}
Let $X$ be a $\mathbf{B}$-valued random variable and let $0 < p <
1$. Then
\begin{equation}
\sum_{n = 1}^{\infty} \frac{1}{n} \left(
\frac{\|S_{n}\|}{n^{1/p}}\right) < \infty~~\mbox{a.s.}
\end{equation}
if and only if
\begin{equation}
\mathbb{E}\|X\|^{p} < \infty.
\end{equation}
Furthermore,
\begin{equation}
\sum_{n = 1}^{\infty} \frac{1}{n} \left(
\frac{\mathbb{E}\|S_{n}\|}{n^{1/p}}\right) < \infty
\end{equation}
if and only if
\begin{equation}
\mathbb{E}\|X\| < \infty.
\end{equation}
\end{theorem}

\vskip 0.3cm

The proofs of Theorems 2.1, 2.2, and 2.4 will be given in Section 4.
Theorem 2.3 and Corollaries 2.2 and 2.3 will be proved
in Section 5. For illustrating the conditions (2.11), (2.12), and
(2.13) of Theorem 2.3, Section 5 also contains some examples.

We now summarize our Theorems 2.1-2.4 and Corollaries 2.1-2.3 for a
real-valued random variable $X$. For $1 < p < 2$, the equivalence of
{\bf (iii)} and {\bf (iv)} has recently been discovered by Hechner
and Heinkel [4] (see Theorem 1.3 above) assuming $\mathbb{E}X = 0$
for the implication ((iii) $\Rightarrow$ (iv)).

\vskip 0.3cm

\begin{theorem}
Let $ \{X_{n};~n \geq 1 \}$ be a sequence of independent copies of a
real-valued random variable $X$. For $0 < p < 2$, the following
two statements are equivalent: \\
\\
$\displaystyle {\bf (i)} \quad \sum_{n = 1}^{\infty}
\frac{1}{n}\left(\frac{|S_{n}|}{n^{1/p}}\right)
< \infty~~\mbox{a.s.}$,\\
\\
$\displaystyle {\bf (ii)} \quad \left \{
\begin{array}{ll}
\mathbb{E}|X|^{p} < \infty, & \mbox{if}~~0 < p < 1\\
&\\
\mathbb{E}X = 0,~\sum_{n=1}^{\infty} \frac{\left|\mathbb{E}XI\{|X|
\leq n\} \right|}{n} < \infty,&\\
&\\
\mbox{and}~\sum_{n=1}^{\infty} \frac{\int_{\min\{u_{n}, n\}}^{n}
\mathbb{P}(|X| > t) dt}{n} <
\infty, & \mbox{if}~~p = 1\\
&\\
\mathbb{E}X = 0 ~~\mbox{and}~~\int_{0}^{\infty}
\mathbb{P}^{1/p}(|X|
> t) dt < \infty, & \mbox{if}~~1 < p < 2.
\end{array}
\right. $\\
\\
For $0 < p < 2$, the following two statements are
equivalent: \\
\\
$\displaystyle {\bf (iii)} \quad \sum_{n = 1}^{\infty} \frac{1}{n}
\left(\frac{\mathbb{E}|S_{n}|}{n^{1/p}}\right)
< \infty$,\\
\\
$\displaystyle {\bf (iv)} \quad \left \{
\begin{array}{ll}
\mathbb{E}|X| < \infty, & \mbox{if}~~0 < p < 1 \\
&\\
\mathbb{E}X = 0 ~~\mbox{and}~~\mathbb{E}|X| \ln(1 + |X|) < \infty,
& \mbox{if}~~p = 1\\
&\\
\mathbb{E}X = 0 ~~\mbox{and}~~\int_{0}^{\infty}
\mathbb{P}^{1/p}(|X|
> t) dt < \infty, & \mbox{if}~~1 < p < 2.
\end{array}
\right. $\\
\\
Furthermore, for $1 < p < 2$, the three statements {\bf (i)}, {\bf
(ii)}, and {\bf (iii)} are equivalent and anyone of them implies
\[
\lim_{n \rightarrow \infty} \frac{S_{n}}{n^{1/p}} = 0 ~~\mbox{a.s.}
\]
Additionally, for the case $p = 1$, if
\[
X~ \mbox{is symmetric with}~ \mathbb{E}|X| \ln^{\delta}(1 + |X|) <
\infty ~~\mbox{for some}~ \delta > 0
\]
or
\[
\mathbb{E}X = 0~~\mbox{and}~~\mathbb{E}|X| \ln(1 + |X|) < \infty,
\]
then
\[
\sum_{n=1}^{\infty} \frac{1}{n} \left(\frac{|S_{n}|}{n} \right) <
\infty ~~\mbox{a.s.}
\]
\end{theorem}

\section{Some general results for
$\sum_{n=1}^{\infty} a_{n}\left\| \sum_{k=1}^{n} V_{k} \right\|$}

In this section, by using the remarkable Hoffmann-J{\o}rgensen [5]
inequality, we establish in Theorem 3.1 some general results for sums of the form
$\sum_{n=1}^{\infty} a_{n}\left\| \sum_{k=1}^{n} V_{k} \right\|
~(a_{n} \in [0, \infty))$. These results will be used for proving
the main results and they may be of independent interest. Theorem 3.1 is a
modified version of the authors' original result and this modification
and its elegant proof were so kindly presented to us by the Referee.

The following lemma will be used in the current section and in Sections 4 and 5.

\vskip 0.3cm

\begin{lemma}
Let $g: ~[0, \infty) \rightarrow [0, \infty)$ be a convex function
and let $\{Y_{1}, ..., Y_{n}\}$ be a set of $n \geq 2$ independent
$\mathbf{B}$-valued random variables such that
$\mathbb{E}g\left(\|Y_{i}\|\right) < \infty$ for $1 \leq i \leq n$.
Then we have:

{\bf (i)}~ If $\mathbb{E}\|Y_{2}\| < \infty$, then
\[
\mathbb{E}g\left(\|Y_{1} - \mathbb{E}Y_{2} \|\right)
\leq \mathbb{E}g \left(\|Y_{1} - Y_{2}\|\right).
\]
Hence, in particular,
\[
\mathbb{E}\|Y_{1} - \mathbb{E}Y_{2} \| \leq \mathbb{E}\|Y_{1} - Y_{2}\|.
\]

{\bf (ii)}~ If $\mathbb{E}Y_{1} = \cdots = \mathbb{E}Y_{n} = 0$, then
\[
\mathbb{E} g\left(\max_{1 \leq k \leq n} \|Y_{k} \| \right)
\leq 2 \mathbb{E}g\left(2\left\|\sum_{k=1}^{n} Y_{k} \right\| \right).
\]
Hence, in particular,
\[
\mathbb{E} \left(\max_{1 \leq k \leq n} \|Y_{k} \| \right) \leq 4
\mathbb{E} \left\|\sum_{k=1}^{n} Y_{k} \right\|.
\]
\end{lemma}

{\it Proof}~ Applying (2.5) of Ledoux and Talagrand [10, p. 46],
part (i) follows immediately. Let $\{R_{1}, ..., R_{n}\}$ be
independent Rademacher random variables independent of $\{Y_{1},
..., Y_{n}\}$. Since $\mathbb{E}Y_{1} = \cdots = \mathbb{E}Y_{n} =
0$, by Proposition 2.3 of Ledoux and Talagrand [10, p. 47] and Lemma
6.3 of Ledoux and Talagrand [10, p. 152], we have
\[
\mathbb{E}g\left(\max_{1 \leq k \leq n} \|Y_{k} \| \right) =
\mathbb{E}g\left(\max_{1 \leq k \leq n} \|R_{k}Y_{k} \| \right) \leq
2 \mathbb{E}g\left(\left\|\sum_{k=1}^{n} R_{k}Y_{k} \right\|\right) \leq 2
\mathbb{E}g\left(2\left\|\sum_{k=1}^{n} Y_{k} \right\| \right)
\]
proving part (ii). ~$\Box$

\vskip 0.3cm

\begin{theorem}
Let $\{V_{n};~n \geq 1 \}$ be a sequence of independent {\bf
B}-valued random variables and let $\{a_{n};~n \geq 1 \}$ be a
sequence of nonnegative numbers such that $\sum_{n=1}^{\infty}a_{n}
< \infty$. Set
\[
b_{n} = \sum_{i=n}^{\infty} a_{i}, ~~T_{n} = \sum_{i=1}^{n} V_{i},
~~n \geq 1, ~~N = \sup_{n \geq 1}b_{n}\left\|V_{n}\right\|, ~~L =
\sum_{n=1}^{\infty}a_{n}\left\|T_{n}\right\|,
\]
\[
M = \sum_{n=1}^{\infty} b_{n} \left\|V_{n}\right\|,
~~\Gamma(t) = \sum_{n=1}^{\infty} \mathbb{P}\left(b_{n}
\left\|V_{n}\right\| > t \right), ~~t \geq 0.
\]
Then we have $L \leq M$. Suppose that $\sup_{n \geq 1}
\mathbb{P}\left(b_{n} \|V_{n}\| > a\right) < 1$ for some $a \geq 0$.
Then we have the following five conclusions:
\begin{align*}
& {\bf (a)} \quad L < \infty~\mbox{a.s.} \Longrightarrow N < \infty
~\mbox{a.s.} \Longleftrightarrow \inf_{r \geq 0} \Gamma(r) <
\infty. \\
& {\bf (b)} \quad \mathbb{E}N < \infty \Longleftrightarrow
\int_{r}^{\infty} \Gamma(t)dt < \infty ~\mbox{for some}~r \geq
0.\\
& {\bf (c)} \quad \mathbb{E}L < \infty \Longleftrightarrow L <
\infty ~\mbox{a.s.}~ \mbox{and}~ \mathbb{E}N <
\infty.\\
& {\bf (d)} \quad \mathbb{E}L < \infty \Longrightarrow
b_{n} \mathbb{E}\left\|V_{n}\right\| < \infty ~\forall ~n \geq 1~
\mbox{and}~ \lim_{n \rightarrow \infty}b_{n} \mathbb{E} \left\|T_{n} -
\mathbb{E}T_{n} \right\| = 0. \\
& {\bf (e)} \quad \sum_{n=1}^{\infty}b_{n} = \infty ~\mbox{and}~ L <
\infty ~\mbox{a.s.} \Longrightarrow \liminf_{n \rightarrow
\infty}\frac{\left\|T_{n}\right\|}{n} = 0~ \mbox{a.s.}
\end{align*}
\end{theorem}

\vskip 0.3cm

{\it Proof}~ Clearly we have
\[
L \leq \sum_{n=1}^{\infty} a_{n}\sum_{i=1}^{n} \left\|V_{i}\right\|
= \sum_{i=1}^{\infty}\left(\sum_{i=n}^{\infty}a_{i}\right)\left\|V_{i}\right\|
= \sum_{i=1}^{\infty}b_{i}\left\|V_{i}\right\| = M.
\]
By the Kolmogorov $0$-$1$ law, we have $\mathbb{P}(N = \infty) = 0$ or $1$ and since
$1 - x \leq e^{-x}$, we have that, for all $t \geq 0$,
\[
1 - e^{-\Gamma(t)} \leq
1 - \prod_{n=1}^{\infty}\mathbb{P}\left(b_{n}\left\|V_{n}\right\| \leq t \right)
= \mathbb{P}(N > t) \leq \Gamma(t).
\]
Hence, we see that the equivalence in {\bf (a)} holds. We will verify below that
$L < \infty$ a.s. $\Longrightarrow N < \infty$ a.s.

Suppose that $\Gamma(r) < \infty$. Since $\Gamma(t)$ and $h(x) =
\left(1 - e^{-x} \right)/x$ are decreasing (with $h(0) :=1$), we
have $\gamma \Gamma(t) \leq 1 - e^{-\Gamma(t)}$ for all $t \geq r$
where $\gamma = h(\Gamma(r))$ and since $\gamma > 0$, we see that
{\bf (b)} holds.

Let $\{V^{'}_{n};~n \geq 1 \}$ be an independent copy of $\{V_{n};~n
\geq 1 \}$ and let $\hat{V}_{n} = V_{n} - V^{'}_{n}, ~n \geq 1$. Let
$\hat{T}_{n}, ~\hat{L}, ~\hat{N}$ and $\hat{\Gamma}(t)$ be defined
as above with $\{V_{n};~n \geq 1 \}$ replaced by $\{\hat{V}_{n};~n
\geq 1 \}$. By assumption there exist $\lambda > 0, a > 0$ such that
$\mathbb{P}\left(b_{n}\|V_{n}\| \leq a \right) \geq \lambda$ for all
$n \geq 1$. Hence, by the symmetrization inequality we have $\lambda
\Gamma(t + a) \leq \hat{\Gamma}(t)$ for all $t \geq 0$; see (6.1) in
[10, p. 150]. Let $\|x\|_{1} := \sum_{i=1}^{n}\|x_{i}\| < \infty$
denote the $\ell^{1}$-norm on $\mathbf{B}^{n}$. Then
$(\mathbf{B}^{n}, \|\cdot\|_{1})$ is a real separable Banach space
and if we set
\[
Y_{i}^{(n)} = \left(0, ..., 0, a_{i}\hat{V}_{i}, ..., a_{n}\hat{V}_{i}\right),~~
i = 1, ..., n,
\]
then $Y_{1}^{(n)}, ..., Y_{n}^{(n)}$ are independent symmetric $\mathbf{B}^{n}$-valued
random variables satisfying
\[
\|Y_{i}^{(n)}\|_{1} = \left(\sum_{k=i}^{n}a_{k}\right)\|\hat{V}_{i}\|
= \left(b_{i} - b_{n+1}\right) \|\hat{V}_{i}\|, ~~i = 1, ..., n,
~~\left\|\sum_{i=1}^{n}Y_{i}^{(n)}\right\|_{1}
= \sum_{i=1}^{n}a_{i} \|\hat{T}_{i}\|.
\]
Hence, by Proposition 2.3 in [10, p. 47] and (3.3) in [6] we have
that, for all $t \geq 0$,
\[
\mathbb{P}\left(\hat{N}_{n} > t \right) \leq 2 \mathbb{P}\left(\hat{L}_{n} > t\right),
~~\mathbb{P}\left(\hat{L}_{n} > 3t \right) \leq \mathbb{P}\left(\hat{N}_{n} > t\right)
+ 4 \mathbb{P}^{2} \left(\hat{L}_{n} > t \right),
\]
where
\[
\hat{N}_{n} = \max_{1 \leq i \leq n}\|Y_{i}^{(n)}\|_{1}
= \max_{1 \leq i \leq n} \left(b_{i} - b_{n+1}\right) \|\hat{V}_{i}\|
~~\mbox{and}~~ \hat{L}_{n} = \left\|\sum_{i=1}^{n}Y_{i}^{(n)}\right\|_{1} =
\sum_{i=1}^{n} a_{i} \|\hat{T}_{i}\|.
\]
Since $\hat{N}_{n} \nearrow \hat{N}$ and $\hat{L}_{n} \nearrow \hat{L}$ as $n \rightarrow \infty$,
we have that, for all $t \geq 0$,
\begin{equation}
\mathbb{P}\left(\hat{N} > t \right) \leq 2 \mathbb{P}\left(\hat{L} > t\right),
~~\mathbb{P}\left(\hat{L} > 3t \right) \leq \mathbb{P}\left(\hat{N} > t\right)
+ 4 \mathbb{P}^{2} \left(\hat{L} > t \right).
\end{equation}
Suppose that $L < \infty$ a.s. We then have $\hat{L} < \infty$ a.s. and so by (3.1) we get
$\hat{N} < \infty$ a.s. Hence, by the Borel-Cantelli lemma, there exists $r \geq 0$ such
that $\hat{\Gamma}(r) < \infty$ and since $\lambda \Gamma(r + a) \leq \hat{\Gamma}(r)$,
we have $N < \infty$ a.s. which proves the first implication in {\bf (a)}.

Suppose that $\mathbb{E}L < \infty$. Then we have $L < \infty$ a.s. and $\mathbb{E}\hat{L} < \infty$
and so by (3.1), we have $\mathbb{E}\hat{N} < \infty$. Hence, by {\bf (b)} there exists $r > 0$ such
that $\int_{r}^{\infty} \hat{\Gamma}(t)dt < \infty$ and since
\[
\lambda \mathbb{P}(N > t + a)
\leq \lambda \Gamma(t + a) \leq \hat{\Gamma}(t),
\]
we have $\mathbb{E}N < \infty$ which proves the
implication ``$\Longrightarrow$" in {\bf (c)}.

Suppose that $L < \infty$ a.s. and $\mathbb{E}N < \infty$.
Then we have $\mathbb{E}\hat{N} < \infty$ and $\hat{L} < \infty$ a.s.
So by (3.1) and a standard argument (see [6]) we have $\mathbb{E}\hat{L} < \infty$.
By Lemma 3.1 (i), we have $\mathbb{E}\left\|T_{n} - \mathbb{E}T_{n}\right\|
\leq \mathbb{E}\left\|\hat{T}_{n}\right\|, ~n \geq 1$ and so we have
\[
\sum_{n=1}^{\infty} a_{n} \mathbb{E}\left\|T_{n} - \mathbb{E}T_{n}\right\|
\leq  \sum_{n=1}^{\infty} a_{n} \mathbb{E}\|\hat{T}_{n}\|
= \mathbb{E}\hat{L} < \infty.
\]
In particular, we have
\[
\sum_{n=1}^{\infty} a_{n} \|T_{n} - \mathbb{E}T_{n}\|
< \infty ~~\mbox{a.s.}
\]
and since $L = \sum_{n=1}^{\infty}a_{n} \|T_{n}\| < \infty$ a.s., we
have $\sum_{n=1}^{\infty} a_{n} \left\|\mathbb{E}T_{n}\right\| <
\infty$. Since $\mathbb{E}\left\|T_{n}\right\| \leq
\mathbb{E}\left\|T_{n} - \mathbb{E}T_{n} \right\| +
\left\|\mathbb{E}T_{n} \right\|,~ n \geq 1$, we have $\mathbb{E}L =
\sum_{n=1}^{\infty} a_{n} \mathbb{E}\left\|T_{n}\right\| < \infty$
which proves the implication ``$\Longleftarrow$" in {\bf (c)}.

Suppose that $\mathbb{E}L < \infty$. Then we have
$\mathbb{E}\left\|V_{n}\right\| < \infty$ and by Lemma 3.1 (i), we have $\mathbb{E}\left\|T_{n} -
\mathbb{E}T_{n} \right\| \leq \mathbb{E}\left\|T_{i} -
\mathbb{E}T_{i} \right\| \leq 2 \mathbb{E}\|T_{i} \|$ for
all $1 \leq n \leq i$. Hence, we have
\[
b_{n} \mathbb{E}\left\|T_{n} - \mathbb{E}T_{n} \right\| =
\sum_{i=n}^{\infty} a_{i}\mathbb{E}\left\|T_{n} - \mathbb{E}T_{n}
\right\| \leq 2 \sum_{i=n}^{\infty} a_{i}
\mathbb{E}\left\|T_{i}\right\|
\]
and since $\sum_{i=1}^{\infty} a_{i} \mathbb{E}\left\|T_{i}\right\|
= \mathbb{E}L < \infty$, we see that $b_{n}\mathbb{E} \left\|T_{n} -
\mathbb{E}T_{n} \right\| \rightarrow 0$ as $n \rightarrow \infty$.

Since $\sum_{n=1}^{\infty} b_{n} = \sum_{n=1}^{\infty} na_{n}$, we
see that {\bf (e)} holds. ~$\Box$

\section{Proofs of Theorems 2.1, 2.2, and 2.4}

Let $\{X_{n};~n \geq 1 \}$ be a sequence of independent copies
of a $\mathbf{B}$-valued random variable $X$. We consider a special case
of Theorem 3.1 which will be used in the proofs of Theorems 2.1-2.4 and
Corollary 5.1. Set $V_{n} = X_{n}$, $T_{n} = S_{n}$, $n \geq 1$. For given
$0 < p < 2$, write
\[
G(t) = \mathbb{P}(\|X\| > t)~~\mbox{for}~t \geq 0~~\mbox{and}~~
a_{n} = n^{-1/p} - (n + 1)^{-1/p}, ~n \geq 1.
\]
Then the hypotheses of Theorem 3.1 hold and by the mean value
theorem we have
\begin{equation}
\left(p^{-1}2^{-1 - 1/p}\right) n^{-1 - 1/p} = p^{-1}(2n)^{-1 - 1/p}
\leq p^{-1}(n+1)^{-1 - 1/p} \leq a_{n} \leq \left(p^{-1}\right)n^{-1
- 1/p}, ~~n \geq 1.
\end{equation}
Let us define
\[
\Lambda_{p}(x) = \sum_{n=0}^{\infty} (n+1)^{-1/p}x^{n}, ~~0 \leq x
\leq 1.
\]
By partial summation, we have
\[
(1-x) \Lambda_{p}(x) = \sum_{n=1}^{\infty} \left(1 - x^{n}\right)
a_{n}
\]
and since
\[
\mathbb{P}\left(\max_{1 \leq i \leq n} \|X_{i} \| > t \right) = 1 -
\left(1 - G(t)\right)^{n},
\]
we have
\begin{equation}
\begin{array}{lll}
\mbox{$\displaystyle \sum_{n=1}^{\infty} a_{n}
\mathbb{E}\left(\max_{1 \leq i \leq n} \left\|X_{i}\right\|\right)$}
&=& \mbox{$\displaystyle \sum_{n=1}^{\infty} a_{n} \int_{0}^{\infty}
\left(1 - (1 - G(t))^{n} \right)dt $}\\
&&\\
&=& \mbox{$\displaystyle \int_{0}^{\infty} \left(\sum_{n=1}^{\infty}
a_{n}
\left(1 - (1 - G(t))^{n} \right)\right)dt $}\\
&&\\
&=& \mbox{$\displaystyle \int_{0}^{\infty} G(t) \Lambda_{p}(1 -
G(t))dt.$}
\end{array}
\end{equation}
Adopting the notation of Theorem 3.1, we have
\[
b_{n} = n^{-1/p}, ~n \geq 1~~\mbox{and}~~\Gamma(t) =
\sum_{n=1}^{\infty}\mathbb{P}\left(\|X\| > t n^{1/p} \right).
\]
Let $r > 0$ be given. Since, for every real-valued random variable
$U$, $\mathbb{E}(U - r)^{+} = \int_{r}^{\infty} \mathbb{P}(U > t)
dt$, we have
\[
\begin{array}{lll}
\mbox{$\displaystyle \int_{r}^{\infty} \Gamma(t)dt$}
&=&
\mbox{$\displaystyle \sum_{n=1}^{\infty} \int_{r}^{\infty}
\mathbb{P}\left(\|X\| > t n^{1/p} \right) dt$}\\
&&\\
&=& \mbox{$\displaystyle \sum_{n=1}^{\infty}
\mathbb{E}\left(\frac{\|X\|}{n^{1/p}} - r \right)^{+}$}\\
&&\\
&=& \mbox{$\displaystyle r \mathbb{E}\left(\sum_{n=1}^{\infty}
\left(\frac{\|X\|}{rn^{1/p}} - 1 \right)^{+}\right)$}\\
&&\\
&=& \mbox{$\displaystyle r \mathbb{E} \phi_{p}
\left(\frac{1}{r}\|X\|\right)$}
\end{array}
\]
where
\[
\phi_{p}(s) = \sum_{n=1}^{\infty} \left(s n^{-1/p} - 1 \right)^{+} =
\sum_{1 \leq n \leq s^{p}} \left(s n^{-1/p} - 1 \right), ~~s \geq 0.
\]
Clearly, $\phi_{p}(s)$ is nondecreasing such that $\phi_{p}(s) = 0$ for all
$0 \leq s \leq 1$ and $\phi_{p}(s) > 0$ for all
$s > 1$. Let
\[
\psi_{p}(s) = \left\{
\begin{array}{ll}
s^{p},~~&\mbox{if}~ 1 < p < 2 \\
s\ln(1+s),~~& \mbox{if}~p = 1 \\
s,~~&\mbox{if}~0 < p < 1.
\end{array}
\right.
\]
It is easy to see that
\[
\lim_{s \rightarrow \infty} \frac{\phi_{p}(s)}{\psi_{p}(s)} = \left\{
\begin{array}{ll}
1/(p-1),~~&\mbox{if}~ 1 < p < 2 \\
1,~~& \mbox{if}~p = 1\\
\sum_{n=1}^{\infty}n^{-1/p},~~&\mbox{if}~0 < p < 1.
\end{array}
\right.
\]
Thus, there exist positive constants $0 < c_{p} < C_{p}$ satisfying
\begin{equation}
c_{p} \psi_{p}(s) \leq \phi_{p}(s) \leq C_{p} \psi_{p}(s) ~~\mbox{for all}~ s \geq 2.
\end{equation}
By Lemma 3.1 (ii), we have for $n \geq 1$,
\[
\begin{array}{lll}
\mbox{$\displaystyle \mathbb{E}\left(\max_{1 \leq i \leq n} \|X_{i}\| \right)$}
&\leq& \mbox{$\displaystyle \|\mathbb{E}X\|
+ \mathbb{E}\left(\max_{1 \leq i \leq n} \|X_{i} - \mathbb{E}X\| \right)$}\\
&&\\
&\leq& \mbox{$\displaystyle \|\mathbb{E}X\| + 4\mathbb{E}\left\|S_{n} - \mathbb{E}S_{n}\right\|$}\\
&&\\
&\leq& \mbox{$\displaystyle \|\mathbb{E}X\| + 8\mathbb{E}\left\|S_{n}\right\|.$}
\end{array}
\]
Hence, by (4.2), (4.3), and Theorems 3.1 and 1.1, we obtain the following conclusions:
\begin{align*}
& {\bf (C1)} \quad \sum_{n=1}^{\infty} a_{n} \|S_{n}\| < \infty
~\mbox{a.s.} \Longrightarrow \sup_{n \geq 1} n^{-1/p} \|X_{n}\| <
\infty ~\mbox{a.s.}
\Longleftrightarrow \mathbb{E}\|X\|^{p} < \infty.\\
& {\bf (C2)} \quad \mathbb{E}\left(\sup_{n \geq 1} n^{-1/p}
\|X_{n}\| \right) < \infty \Longleftrightarrow
\mathbb{E}\psi_{p}(\|X\|) < \infty \Longrightarrow \mathbb{E}\|X\| < \infty. \\
& {\bf (C3)} \quad \sum_{n=1}^{\infty} a_{n} \mathbb{E}\|S_{n}\| <
\infty \Longleftrightarrow \sum_{n=1}^{\infty} a_{n} \|S_{n}\| <
\infty
~\mbox{a.s.~and}~ \mathbb{E}\psi_{p}(\|X\|) < \infty.\\
& {\bf (C4)} \quad 1 \leq p < 2~~\mbox{and}~~\sum_{n=1}^{\infty}
a_{n} \|S_{n}\| < \infty ~\mbox{a.s.}
\Longrightarrow \mathbb{E}\|X\| < \infty~~\mbox{and}~~\mathbb{E}X = 0. \\
& {\bf (C5)} \quad \mbox{If}~\sum_{n=1}^{\infty} a_{n}
\mathbb{E}\|S_{n}\| < \infty, \mbox{~then~} \int_{0}^{\infty}
G(t)\Lambda_{p}(1 - G(t)) dt < \infty \mbox{~and~}
\frac{S_{n}}{n^{1/p}} \rightarrow 0 ~\mbox{a.s. and in~}
L^{1}(\mathbb{P}).
\end{align*}

\vskip 0.3cm

{\it Proof of Theorem 2.1}~ By (4.1), {\bf (C1)}, and {\bf (C3)}, we
see that (2.1) and (2.2) are equivalent. For the given $1 < p < 2$,
by (2.5) in the proof of Lemma 2 in Hechner and Heinkel [4], there
exists a constant $0 < A_{p} < 1$ such that
\[
A_{p}x^{1/p -1} \leq \Lambda_{p}(1-x) \leq x^{1/p - 1}~~\mbox{for
all}~ 0 < x \leq 1.
\]
Hence, we see that the last part of Theorem 2.1 follows from {\bf
(C1)}, {\bf (C4)}, and {\bf (C5)}. ~$\Box$

\vskip 0.3cm

{\it Proof of Theorem 2.2}~ The theorem is an immediate consequence of {\bf (C3)}
and {\bf (C4)}. ~$\Box$

\vskip 0.3cm

{\it Proof of Theorem 2.4}~ First, by (4.1) and {\bf (C1)}, (2.17) follows from (2.16).

Conversely, suppose that (2.17) holds. By (4.1) and Theorem 3.1, we
have
\[
\left(p^{-1}2^{-1 - 1/p}\right) \sum_{n=1}^{\infty}\frac{1}{n}\left(\frac{\|S_{n}\|}{n^{1/p}}\right)
\leq \sum_{n=1}^{\infty}a_{n} \|S_{n}\| \leq \sum_{n=1}^{\infty}n^{-1/p} \|X_{n}\|.
\]
Since $0 < p < 1$, by (2.17) we have
\[
\sum_{n=1}^{\infty}n^{-1/p} \|X_{n}\| < \infty~~\mbox{a.s.}
\]
(see Theorem 5.1.3 in Chow and Teicher [2, p. 115]) and (2.16)
follows.

Furthermore, for $0 < p < 1$, the equivalence of (2.18) and (2.19) is trivial since
\[
\mathbb{E}\|X\| \leq \sum_{n=1}^{\infty}\frac{1}{n} \left(\frac{\mathbb{E}\|S_{n}\|}{n^{1/p}} \right)
\leq \left(\sum_{n=1}^{\infty} \frac{1}{n^{1/p}} \right) \mathbb{E}\|X\| < \infty.
\]
The proof of Theorem 2.4 is complete. ~$\Box$

\section{Proofs of Theorem 2.3 and Corollaries 2.2 and 2.3}

Let $\{X_{n}; ~n \geq 1 \}$ be a sequence of independent copies of
the $\mathbf{B}$-valued random variable $X$. Write
\[
X_{n}^{(1)} = X_{n} I\{\|X_{n}\| \leq n\}~~\mbox{and}~~X_{n}^{(2)} =
X_{n} I\{\|X_{n}\| > n\}, ~~n \geq 1
\]
and
\[
S_{n}^{(1)} = \sum_{k=1}^{n} X_{k}^{(1)} ~~\mbox{and}~~S_{n}^{(2)} =
\sum_{k=1}^{n} X_{k}^{(2)}, ~~n \geq 1.
\]

For the proof of Theorem 2.3, we need the following four preliminary
lemmas.

\vskip 0.3cm

\begin{lemma}
Let $X$ be a {\bf B}-valued random variable with $\mathbb{E}\|X\| <
\infty$. Then
\begin{equation}
\sum_{n=1}^{\infty} \frac{1}{n^{2}}
\sum_{k=1}^{n}\mathbb{E}\|X\|I\{k < \|X\| \leq n\} < \infty
\end{equation}
and
\begin{equation}
\sum_{n=1}^{\infty} \frac{u_{n}}{n^{2}} < \infty.
\end{equation}
\end{lemma}

{\it Proof}~ Since $\mathbb{E}\|X\| < \infty$, we have that
\[
\begin{array}{lll}
\mbox{$\displaystyle \sum_{n=1}^{\infty} \frac{1}{n^{2}}
\sum_{k=1}^{n}\mathbb{E}\|X\|I\{k < \|X\| \leq n\}$} &=&
\mbox{$\displaystyle \sum_{k=1}^{\infty} \sum_{n=k}^{\infty}
\frac{1}{n^{2}} \mathbb{E}\|X\|I\{k < \|X\| \leq
n\}$}\\
&&\\
&\leq& \mbox{$\displaystyle \sum_{k=1}^{\infty} \sum_{n=k}^{\infty}
\frac{1}{n^{2}} \sum_{m=k}^{n} m \mathbb{P}(m-1 <
\|X\| \leq m)$}\\
&&\\
&=& \mbox{$\displaystyle \sum_{k=1}^{\infty}
\sum_{m=k}^{\infty}m\left(\sum_{n=m}^{\infty} \frac{1}{n^{2}}
\right)
\mathbb{P}(m-1 < \|X\| \leq m)$}\\
&&\\
&\leq& \mbox{$\displaystyle \sum_{k=1}^{\infty} \sum_{m=k}^{\infty}
m \times \frac{2}{m}\mathbb{P}(m-1 < \|X\| \leq m)$}\\
&&\\
&=& \mbox{$\displaystyle 2\sum_{m=1}^{\infty} \sum_{k=1}^{m}
\mathbb{P}(m-1 < \|X\| \leq m)$}\\
&&\\
&=& \mbox{$\displaystyle 2\sum_{m=1}^{\infty} m \mathbb{P}(m-1 < \|X\| \leq m)$}\\
&&\\
&<& \mbox{$\displaystyle \infty$}
\end{array}
\]
proving (5.1).

We now show that (5.2) follows from $\mathbb{E}\|X\| < \infty$.
Since $\{u_{n}; ~n \geq 1 \}$ is an increasing sequence with
\[
\sup_{n \geq 1} u_{n} = \sup\{x:~\mathbb{P}(|X\| \leq x) < 1\},
\]
we obtain that
\[
\begin{array}{lll}
\mbox{$\displaystyle \sum_{n=1}^{\infty} \frac{u_{n+1} - u_{n}}{n}$}
&\leq& \mbox{$\displaystyle \sum_{n=1}^{\infty}
2\int_{u_{n}}^{u_{n+1}} \mathbb{P}(\|X\| > x)dx$}\\
&&\\
&=& \mbox{$\displaystyle 2\int_{u_{1}}^{\infty} \mathbb{P}(\|X\| >
x) dx$} \\
&&\\
&\leq& \mbox{$\displaystyle 2 \mathbb{E}\|X\| < \infty.$}
\end{array}
\]
Note that for each integer $m \geq 2$, partial summation yields
\[
\begin{array}{lll}
\mbox{$\displaystyle \sum_{n=1}^{m} \frac{u_{n+1} - u_{n}}{n}$} &=&
\mbox{$\displaystyle \left(\frac{u_{m+1}}{m} - u_{1} \right) +
\sum_{n=2}^{m} \left(\frac{1}{n-1} - \frac{1}{n}\right) u_{n}$}\\
&&\\
&\geq& \mbox{$\displaystyle -2u_{1} + \sum_{n=1}^{m}
\frac{u_{n}}{n^{2}}.$}
\end{array}
\]
We thus see that
\[
\sum_{n=1}^{m} \frac{u_{n}}{n^{2}} \leq 2u_{1} + \sum_{n=1}^{\infty}
\frac{u_{n+1} - u_{n}}{n} < \infty
\]
proving (5.2). ~$\Box$

\vskip 0.3cm

\begin{lemma}
Let $\{X_{n};~n \geq 1 \}$ be be a sequence of independent copies of
{\bf B}-valued random variable $X$ with $\mathbb{E}\|X\| < \infty$.
Write
\[
U_{n} = \sum_{k=1}^{n} X_{k}I\{\|X_{k}\| \leq n\}, ~n \geq 1.
\]
Then
\[
\begin{array}{ll}
& {\bf (i)} \quad \mbox{$\displaystyle
\sum_{n=1}^{\infty}\frac{1}{n} \left(\frac{\|\mathbb{E}S_{n}^{(1)}
\|}{n} \right) <
\infty~~\mbox{if and only if}~~(2.12)~\mbox{holds},$}\\
&\\
& {\bf (ii)} \quad \mbox{$\displaystyle
\sum_{n=1}^{\infty}\frac{1}{n}
\left(\frac{\mathbb{E}\left\|S_{n}^{(1)} - \mathbb{E}S_{n}^{(1)}
\right\|}{n} \right) < \infty~~\mbox{if and only if}~~
\sum_{n=1}^{\infty}\frac{1}{n} \left(\frac{\mathbb{E}\left\|U_{n} -
\mathbb{E}U_{n} \right\|}{n} \right) < \infty.$}
\end{array}
\]
\end{lemma}

{\it Proof}~ Note that by (5.1) of Lemma 5.1,
\[
\begin{array}{ll}
& \mbox{$\displaystyle
\sum_{n=1}^{\infty}\frac{1}{n}\left\|\frac{\mathbb{E}S_{n}^{(1)}}{n}
- \mathbb{E}XI\{\|X\| \leq n\} \right\|$}\\
&\\
& \mbox{$\displaystyle =
\sum_{n=1}^{\infty}\frac{1}{n}\left\|\frac{\sum_{k=1}^{n}
\left(\mathbb{E}XI\{\|X\| \leq k\} -
\mathbb{E}XI\{\|X\| \leq n\} \right)}{n} \right\|$}\\
&\\
& \mbox{$\displaystyle \leq \sum_{n=1}^{\infty}\frac{1}{n^{2}}
\sum_{k=1}^{n}\mathbb{E}\|X\|I\{k < \|X\| \leq n\}$}\\
&\\
& \mbox{$\displaystyle < \infty.$}
\end{array}
\]
Thus part {\bf (i)} follows.

Similarly, by (5.1) of Lemma 5.1,
\[
\begin{array}{ll}
& \mbox{$\displaystyle \sum_{n=1}^{\infty}\frac{1}{n}
\left(\frac{\mathbb{E}\left\|\left(S_{n}^{(1)} -
\mathbb{E}S_{n}^{(1)}\right) - \left(U_{n} -
\mathbb{E}U_{n}\right)\right\|}{n} \right)$} \\
&\\
& \mbox{$\displaystyle \leq 2\sum_{n=1}^{\infty}\frac{1}{n^{2}}
\sum_{k=1}^{n}\mathbb{E}\|X\|I\{k < \|X\| \leq n\}$} \\
&\\
& \mbox{$\displaystyle < \infty$}
\end{array}
\]
and part {\bf (ii)} follows. ~$\Box$

\vskip 0.3cm

The proof of the next lemma is similar to that of Lemma 4 of Hechner
and Heinkel [4] and contains a nice application of Lemma 1 of
Hechner and Heinkel [4].

\vskip 0.3cm

\begin{lemma}
Let {\bf B} be a Banach space of stable type $1$. Let $\{X_{n};~n
\geq 1 \}$ be a sequence of independent copies of a
$\mathbf{B}$-valued random variable $X$ with $\mathbb{E}\|X\| <
\infty$. Write
\[
U_{n}^{(1)} = \sum_{k=1}^{n} X_{k}I\left\{\|X_{k}\| \leq u_{n}
\right\}, ~n \geq 1.
\]
Then
\begin{equation}
\sum_{n=1}^{\infty}\frac{1}{n}
\left(\frac{\mathbb{E}\left\|U_{n}^{(1)} - \mathbb{E}U_{n}^{(1)}
\right\|}{n} \right) < \infty.
\end{equation}
\end{lemma}

\vskip 0.3cm

{\it Proof}~ Since $\mathbf{B}$ is of stable type $1$, the
Maurey-Pisier [12] theorem asserts that it is also of stable type
$q$ for some $1 < q < 2$. Applying Lemma 1 of Hechner and Heinkel
[4], there exists a universal constant $0 < c(q) < \infty$ such that
\[
\begin{array}{lll}
\mbox{$\displaystyle \mathbb{E}\left\|U_{n}^{(1)} -
\mathbb{E}U_{n}^{(1)} \right\|$} &\leq& \mbox{$\displaystyle c(q)
\left(\sup_{t>0} t^{q} \sum_{k=1}^{n}
\mathbb{P}\left(\|X_{k}\|I\left\{\|X_{k}\| \leq u_{n} \right\} > t
\right) \right)^{1/q}$}\\
&&\\
&\leq& \mbox{$\displaystyle c(q) \left(n \sup_{0 \leq t \leq u_{n}}
t^{q} \mathbb{P}(\|X\| > t) \right)^{1/q}, ~~n \geq 1$}.
\end{array}
\]
It is easy to see that for all $x > 0$,
\[
\begin{array}{lll}
\mbox{$\displaystyle \left(\int_{0}^{x} \mathbb{P}^{1/q}(\|X\| > t)
dt \right)^{q}$} &\geq& \mbox{$\displaystyle \left(\int_{0}^{x}
\mathbb{P}^{1/q}(\|X\|
> x) dt \right)^{q}$}\\
&&\\
&=& \mbox{$\displaystyle x^{q} \mathbb{P}(\|X\|
> x).$}
\end{array}
\]
We thus have that
\[
\begin{array}{lll}
\mbox{$\displaystyle \mathbb{E}\left\|U_{n}^{(1)} -
\mathbb{E}U_{n}^{(1)} \right\|$} &\leq& \mbox{$\displaystyle c(q)
\left(n \sup_{0 \leq t \leq u_{n}} t^{q} \mathbb{P}(\|X\| > t)
\right)^{1/q}$} \\
&&\\
&\leq& \mbox{$\displaystyle c(q) n^{1/q} \int_{0}^{u_{n}}
\mathbb{P}^{1/q}(\|X\|
> t) dt, ~~n \geq 1.$}
\end{array}
\]
Let $u_{0} = 0$ and note that $\mathbb{P}(\|X\| > t) \geq 1/k$ for
$t \in [u_{k-1}, u_{k}), ~k \geq 1$. It follows that
\[
\begin{array}{lll}
\mbox{$\displaystyle \sum_{n=1}^{\infty}\frac{1}{n}
\left(\frac{\mathbb{E}\left\|U_{n}^{(1)} - \mathbb{E}U_{n}^{(1)}
\right\|}{n} \right)$} &\leq& \mbox{$\displaystyle c(q)
\sum_{n=1}^{\infty} \frac{1}{n^{2-1/q}} \int_{0}^{u_{n}}
\mathbb{P}^{1/q}(\|X\|
> t) dt$}\\
&&\\
&=& \mbox{$\displaystyle c(q)\sum_{n=1}^{\infty} \frac{1}{n^{2-1/q}}
\sum_{k=1}^{n}\int_{u_{k-1}}^{u_{k}} \mathbb{P}^{1/q}(\|X\|
> t) dt$}\\
&&\\
&=& \mbox{$\displaystyle c(q)\sum_{k=1}^{\infty}
\left(\sum_{n=k}^{\infty} \frac{1}{n^{2-1/q}} \right)
\int_{u_{k-1}}^{u_{k}} \mathbb{P}^{1/q}(\|X\|
> t) dt$}\\
&&\\
&\leq& \mbox{$\displaystyle C \sum_{k=1}^{\infty}\frac{1}{k^{1-1/q}}
\int_{u_{k-1}}^{u_{k}} \mathbb{P}^{1/q}(\|X\| > t) dt$}\\
&\leq& \mbox{$\displaystyle C \int_{0}^{\infty}\mathbb{P}(\|X\| >
t) dt$} \\
&&\\
&=& \mbox{$\displaystyle C\mathbb{E}\|X\| < \infty$}
\end{array}
\]
proving (5.3) and completing the proof of Lemma 5.3.~$\Box$

\vskip 0.3cm

\begin{lemma}
Let $Y_{1}, Y_{2}, ..., Y_{n}$ be i.i.d. nonnegative real-valued
random variables such that
\begin{equation}
\mathbb{P}\left(Y_{1} > 0 \right) \leq \frac{1}{n}.
\end{equation}
Then
\begin{equation}
\mathbb{E}\left(\max_{1 \leq k \leq n} Y_{k} \right) \geq
\frac{n}{2} \mathbb{E}Y_{1}.
\end{equation}
\end{lemma}

{\it Proof}~ Since $Y_{1}, Y_{2}, ..., Y_{n}$ are i.i.d. nonnegative real-valued
random variables, we have
\[
\mathbb{P} \left(\max_{1 \leq k \leq n}Y_{k} > t \right)
= 1 - \left(1 - \mathbb{P}\left(Y_{1} > t \right)\right)^{n}
\geq 1 - e^{-n \mathbb{P}\left(Y_{1} > t \right)}~~\mbox{for all}~ t \geq 0.
\]
Since $(1 - e^{-x})/x$ is decreasing, we have $1 - e^{-x} \geq \beta x$ for all
$0 \leq x \leq 1$ where $\beta = 1 - e^{-1} \geq 1/2$. Hence, it follows from (5.4)
that
\[
\mathbb{P}\left(\max_{1 \leq k \leq n}Y_{k} > t \right)
\geq \frac{n}{2} \mathbb{P}\left(Y_{1} > t \right)~~\mbox{for all}~ t > 0,
\]
which ensures (5.5). ~$\Box$

\vskip 0.3cm

{\it Proof of Theorem 2.3} ~({\bf Sufficiency}) Since
\[
\|S_{n}\| \leq \|\mathbb{E}S_{n}^{(1)}\| + \|S_{n}^{(1)} -
\mathbb{E}S_{n}^{(1)}\| + \|S_{n}^{(2)}\|, ~n \geq 1,
\]
(2.10) will follow if we can show that (2.11), (2.12), and (2.13)
imply
\begin{equation}
\sum_{n=1}^{\infty} \frac{1}{n}
\left(\frac{\|\mathbb{E}S_{n}^{(1)}\|}{n} \right) < \infty,
\end{equation}
\begin{equation}
\sum_{n=1}^{\infty} \frac{1}{n} \left(\frac{\mathbb{E}\|S_{n}^{(1)}
- \mathbb{E}S_{n}^{(1)}\|}{n} \right) < \infty,
\end{equation}
and
\begin{equation}
\sum_{n=1}^{\infty} \frac{1}{n} \left(\frac{\|S_{n}^{(2)} \|}{n}
\right) < \infty ~~\mbox{a.s.}
\end{equation}
By Lemma 5.2 (i), (5.6) follows from (2.11) and (2.12). Since
\[
\sum_{n=1}^{\infty} \mathbb{P} \left(X_{n}^{(2)} \neq 0 \right) =
\sum_{n=1}^{\infty} \mathbb{P} (\|X\| > n) \leq \mathbb{E}\|X\| <
\infty,
\]
by the Borel-Cantelli lemma we have that
\[
\mathbb{P} \left(X_{n}^{(2)} \neq 0 ~\mbox{i.o.}(n) \right) = 0
\]
which ensures that
\[
\|S_{n}^{(2)}\| = {\it O}\left(1\right) ~~\mbox{a.s. as}~ n
\rightarrow \infty.
\]
Thus (5.8) holds.

We now show that (2.11) and (2.13) imply (5.7). Since
$\mathbb{E}\|X\| < \infty$, by Lemma 5.2 (ii), (5.7) is equivalent
to
\begin{equation}
\sum_{n=1}^{\infty} \frac{1}{n} \left(\frac{\mathbb{E}\|U_{n} -
\mathbb{E}U_{n}\|}{n} \right) < \infty.
\end{equation}
Now (2.9) holds recalling the implication ((2.8) $\Rightarrow$
(2.9)). Hence we can assume, without loss of generality, that $u_{n}
< n$ for all $n \geq 1$. Write
\[
U_{n}^{(2)} = \sum_{k=1}^{n} X_{k}I\left\{u_{n} < \|X_{k}\| \leq n
\right\}, ~n \geq 1.
\]
Clearly, (5.9) will follow provided we can show
\begin{equation}
\sum_{n=1}^{\infty} \frac{1}{n} \left(\frac{\mathbb{E}\|U_{n}^{(1)}
- \mathbb{E}U_{n}^{(1)}\|}{n} \right) < \infty
\end{equation}
and
\begin{equation}
\sum_{n=1}^{\infty} \frac{1}{n} \left(\frac{\mathbb{E}\|U_{n}^{(2)}
- \mathbb{E}U_{n}^{(2)}\|}{n} \right) < \infty.
\end{equation}
Since $\mathbb{E}\|X\| < \infty$ and {\bf B} is of stable type $1$,
by Lemma 5.3, (5.10) holds. Note that, for all $n \geq 1$,
\[
\begin{array}{lll}
\mbox{$\displaystyle \mathbb{E}\|U_{n}^{(2)} -
\mathbb{E}U_{n}^{(2)}\|$} &\leq& \mbox{$\displaystyle
\sum_{k=1}^{n}\mathbb{E}\left\|X_{k}I\left\{u_{n} < \|X_{k}\| \leq n
\right\} - \mathbb{E}X_{k}I\left\{u_{n} < \|X_{k}\| \leq n \right\}
\right\|$}\\
&&\\
&\leq& \mbox{$\displaystyle 2n \mathbb{E}\|X\|I\{u_{n} < \|X\| \leq
n\}$}\\
&&\\
&=& \mbox{$\displaystyle 2n\int_{u_{n}}^{n} t d\mathbb{P}(\|X\| \leq
t)$}\\
&&\\
&\leq& \mbox{$\displaystyle 2nu_{n} \mathbb{P}\left(\|X\| >
u_{n}\right) + 2n\int_{u_{n}}^{n} \mathbb{P}(\|X\| > t) dt$}\\
&&\\
&\leq& \mbox{$\displaystyle 2u_{n} + 2n\int_{u_{n}}^{n}
\mathbb{P}(\|X\| > t) dt.$}
\end{array}
\]
Now (5.2) holds by Lemma 5.1. Thus (5.11) follows from (5.2) and
(2.13). The proof of the sufficiency half of Theorem 2.3 is
complete. ~$\Box$

\vskip 0.3cm

{\it Proof of Theorem 2.3} ~({\bf Necessity}) First, by (4.1) with $p = 1$,
(2.10) is equivalent to
\[
\sum_{n=1}^{\infty} a_{n} \left\|S_{n} \right\| < \infty~~\mbox{a.s.}
\]
Thus, by {\bf (C4)} with $p = 1$, (2.11) follows. In particular, we have
$\mathbb{E}\|X\| < \infty$. Again by (4.1) with $p = 1$ and arguing as in
the proof of (5.8) in the sufficiency half, we have
\[
\sum_{n=1}^{\infty} a_{n} \|S_{n}^{(2)} \| < \infty~~\mbox{a.s.}
\]
Note that
\[
\|S_{n}^{(1)}\| \leq \|S_{n}\| + \|S_{n}^{(2)} \|, ~~n \geq 1.
\]
It thus follows from (2.10) that
\begin{equation}
\sum_{n=1}^{\infty} a_{n} \|S_{n}^{(1)} \| < \infty~~\mbox{a.s.}
\end{equation}
Since, with $p = 1$,
\[
N = \sup_{n \geq 1} b_{n} \|X_{n}^{(1)} \|
= \sup_{n \geq 1} n^{-1} \|X_{n}^{(1)} \| \leq 1,
\]
it follows from (5.12) and Theorem 3.1 (c) that
\[
\sum_{n=1}^{\infty} a_{n} \mathbb{E} \|S_{n}^{(1)} \| < \infty
\]
which by (4.1) ensures that
\begin{equation}
\sum_{n=1}^{\infty} \frac{1}{n}
\left(\frac{\mathbb{E}\|S_{n}^{(1)}\|}{n} \right) < \infty.
\end{equation}
Note that (5.13) implies
\[
\sum_{n=1}^{\infty} \frac{1}{n}
\left(\frac{\|\mathbb{E}S_{n}^{(1)}\|}{n} \right) < \infty
\]
which by $\mathbb{E}\|X\| < \infty$ and Lemma 5.2 (i) yields (2.12).
Moreover, by Lemma 5.2 (ii), (5.13) then ensures that
\begin{equation}
\sum_{n=1}^{\infty} \frac{1}{n} \left(\frac{\mathbb{E}\|U_{n} -
\mathbb{E}U_{n}\|}{n} \right) < \infty.
\end{equation}
Since $\|U_{n}^{(2)} - \mathbb{E}U_{n}^{(2)}\| \leq \|U_{n} -
\mathbb{E}U_{n}\| + \|U_{n}^{(1)} - \mathbb{E}U_{n}^{(1)}\|, ~n \geq
1$ and {\bf B} is of stable type $1$, it follows from Lemma 5.3 and
(5.14) that
\begin{equation}
\sum_{n=1}^{\infty} \frac{1}{n} \left(\frac{\mathbb{E}\|U_{n}^{(2)}
- \mathbb{E}U_{n}^{(2)}\|}{n} \right) < \infty.
\end{equation}
By Lemma 3.1 {\bf (ii)},
\[
\mathbb{E} \max_{1 \leq k \leq n} \left\| X_{k}I\{u_{n} < \|X_{k}\|
\leq n\} - \mathbb{E}XI\{u_{n} < \|X\| \leq n\} \right\| \leq 4
\mathbb{E}\|U_{n}^{(2)} - \mathbb{E}U_{n}^{(2)} \|, ~n \geq 1.
\]
It thus follows from (5.15) and $\mathbb{E}\|X\| < \infty$ that
\[
\begin{array}{ll}
& \mbox{$\displaystyle \sum_{n=1}^{\infty}
\frac{1}{n}\left(\frac{\mathbb{E} \max_{1 \leq k \leq n}
\left\|X_{k}I\{u_{n} < \|X_{k}\| \leq n\}\right\|}{n} \right)$}\\
&\\
& \mbox{$\displaystyle \leq \sum_{n=1}^{\infty}
\frac{1}{n}\left(\frac{4 \mathbb{E}\|U_{n}^{(2)} -
\mathbb{E}U_{n}^{(2)} \|}{n} \right) +
\sum_{n=1}^{\infty}\frac{\mathbb{E}\|X\|}{n^{2}}$}\\
&\\
& \mbox{$\displaystyle < \infty,$}
\end{array}
\]
and hence, by Lemma 5.4, noting that $\mathbb{P}(\|X\| > u_{n}) \leq
n^{-1}, ~n \geq 1$, we get that
\begin{equation}
\sum_{n=1}^{\infty} \frac{\mathbb{E}\|X\|I\{u_{n} < \|X\| \leq
n\}}{n} < \infty.
\end{equation}
Using partial integration, one can easily see that
\begin{equation}
\begin{array}{ll}
& \mbox{$\displaystyle \left|\mathbb{E}\|X\|I\{u_{n} < \|X\| \leq
n\} - \int_{u_{n}}^{n} \mathbb{P}(\|X\| > t)dt \right|$}\\
&\\
& \mbox{$\displaystyle \leq \frac{u_{n}}{n} + n\mathbb{P}(\|X\| >
n), ~~ n \geq 1.$}
\end{array}
\end{equation}
Since $\mathbb{E}\|X\| < \infty$, we have
\begin{equation}
\sum_{n=1}^{\infty} \frac{n\mathbb{P}(\|X\| > n)}{n} =
\sum_{n=1}^{n} \mathbb{P}(\|X\| > n) < \infty,
\end{equation}
and, by Lemma 5.1, (5.2) holds. We thus see that (2.13) follows from
(5.16), (5.17), (5.18), and (5.2) thereby completing the proof of
the necessity half of Theorem 2.3.~ $\Box$

\vskip 0.3cm

\begin{lemma}
Let $X$ be a {\bf B}-valued random variable such that
\[
\mathbb{E} X = 0 ~~\mbox{and}~~\mathbb{E}\|X\|\ln(1 + \|X\|) <
\infty.
\]
Then (2.12) holds.
\end{lemma}

\noindent {\it Proof} ~Argue as in Exercise 5.1.6 (ii) of Chow and
Teicher [2, p. 123]. The details are left to the reader. ~$\Box$

\vskip 0.3cm

\begin{lemma}
Let $X$ be a {\bf B}-valued random variable such that
\begin{equation}
\mathbb{E}\|X\|\ln^{\delta}(1 + \|X\|) < \infty~~\mbox{for
some}~\delta > 0.
\end{equation}
Then (2.13) holds.
\end{lemma}

{\it Proof} ~It follows from (5.17) that, under the condition
$\mathbb{E}\|X\| < \infty$, (2.13) and (5.16) are equivalent. We
thus only need to show that (5.19) implies (5.16). To this end, let
$q = (2+\delta)/2$. Applying H\"{o}lder's inequality, we have that
\[
\begin{array}{ll}
& \mbox{$\displaystyle \mathbb{E}\|X\|I\{u_{n} < \|X\| \leq n\}$}\\
&\\
& \mbox{$\displaystyle = \mathbb{E}\left(\|X\|I\{u_{n} < \|X\|
\leq
n\}\right) \left(I\{u_{n} < \|X\| \leq n\}\right)$}\\
&\\
& \mbox{$\displaystyle \leq \left(\mathbb{E}\left(\|X\|I\{u_{n} <
\|X\| \leq n\}\right)^{q} \right)^{1/q}
\left(\mathbb{E}\left(I\{u_{n} < \|X\| \leq n\}\right)^{q/(q-1)}
\right)^{(q-1)/q} $}\\
&\\
& \mbox{$\displaystyle \leq \left(\mathbb{E}\left(\|X\|I\{\|X\| \leq
n\}\right)^{q} \right)^{1/q} \left(\mathbb{E}I\{\|X\| > u_{n}\}\right)^{(q-1)/q} $}\\
&\\
& \mbox{$\displaystyle = \left(\mathbb{E}\left(\|X\|I\{\|X\| \leq
n\}\right)^{q} \right)^{1/q} \mathbb{P}^{(q-1)/q}\left(\|X\| > u_{n}
\right)$}\\
&\\
& \mbox{$\displaystyle \leq
\frac{\left(\mathbb{E}\left(\|X\|I\{\|X\| \leq n\}\right)^{q}
\right)^{1/q}}{n^{1 - 1/q}}, ~~n \geq 1.$}
\end{array}
\]
Thus, (5.16) follows if we can show that (5.19) implies
\begin{equation}
\sum_{n=1}^{\infty} \frac{\left(\mathbb{E}\left(\|X\|I\{\|X\| \leq
n\}\right)^{q} \right)^{1/q}}{n^{2 - 1/q}} < \infty.
\end{equation}
Write
\[
f_{n} = \frac{\left(\ln^{\delta/q}(1 +
n)\right)\left(\mathbb{E}\|X\|^{q}I\{\|X\| \leq n \}
\right)^{1/q}}{n}, ~~g_{n} = \frac{1}{n^{1-
1/q}\ln^{\delta/q}(1+n)}, ~n \geq 1.
\]
Applying H\"{o}lder's inequality again, we have that
\[
\sum_{n=1}^{\infty} \frac{\left(\mathbb{E}\|X\|^{q}I\{\|X\| \leq n
\} \right)^{1/q}}{n^{2 - 1/q}} = \sum_{n=1}^{\infty} f_{n}g_{n} \leq
\left(\sum_{n=1}^{\infty} f_{n}^{q}
\right)^{1/q}\left(\sum_{n=1}^{\infty}
g_{n}^{q/(q-1)}\right)^{(q-1)/q}.
\]
It is easy to see that
\[
\sum_{n=1}^{\infty} g_{n}^{q/(q-1)} = \sum_{n=1}^{\infty}\frac{1}{n
\ln^{2}(1 + n)} < \infty.
\]
Since $\mathbb{E}\|X\|\ln^{\delta}(1 + \|X\|) < \infty$, we have
that
\[
\begin{array}{lll}
\mbox{$\displaystyle \sum_{n=1}^{\infty} f_{n}^{q}$} &=&
\mbox{$\displaystyle
\sum_{n=1}^{\infty}\frac{(\ln^{\delta}(1+n))\mathbb{E}\|X\|^{q}I\{\|X\|
\leq
n\}}{n^{q}}$}\\
&&\\
&\leq& \mbox{$\displaystyle \sum_{n=1}^{\infty}
\frac{\ln^{\delta}(1+n)}{n^{q}} \sum_{k=1}^{n} k^{q}\mathbb{P}(k-1 <
\|X\|
\leq k)$}\\
&&\\
&=&\mbox{$\displaystyle \sum_{k=1}^{\infty}\left(\sum_{n=k}^{\infty}
\frac{\ln^{\delta}(1+n)}{n^{q}}\right)k^{q}\mathbb{P}(k-1 < \|X\|
\leq k)$}\\
&&\\
&\leq& \mbox{$\displaystyle C\sum_{k=1}^{\infty}
\frac{\ln^{\delta}(1+k)}{k^{q-1}} \times k^{q} \mathbb{P}(k-1 <
\|X\|
\leq k)$}\\
&&\\
&=& \mbox{$\displaystyle C \sum_{k=1}^{\infty}(k
\ln^{\delta}(1+k))\mathbb{P}(k-1 < \|X\|
\leq k)$}\\
&&\\
&<& \mbox{$\displaystyle \infty.$}
\end{array}
\]
Thus (5.20) holds. This completes the proof of Lemma 5.6. ~$\Box$

\vskip 0.3cm

{\it Proof of Corollary 2.2} ~By Theorem 2.2, we only need to show
that, under the assumption that $\mathbf{B}$ is of stable type $1$,
(2.6) follows from (2.7). Clearly, (2.7) implies (2.11). Applying
Lemmas 5.5 and 5.6, (2.12) and (2.13) follow from the second half of
(2.7). Thus, by Theorem 2.3, (2.6) holds. ~$\Box$

\vskip 0.3cm

{\it Proof of Corollary 2.3} ~The conclusion of Corollary 2.3
follows immediately from Lemma 5.6 and Theorem 2.3. ~$\Box$

\vskip 0.3cm

From the definition of $\Lambda_{1}(x)$ given in Section 4, we see that
\[
\Lambda_{1}(0) = 1 ~~\mbox{and}~~\Lambda_{1}(x) = - \frac{\ln(1-x)}{x} ~~
\mbox{for all}~ 0 < x \leq 1.
\]
Hence, by {\bf (C5)} we see that (2.5) implies
\begin{equation}
\int_{0}^{\infty}G(t) \ln \frac{1}{G(t)} dt < \infty.
\end{equation}
Conversely, if (5.21) holds, then (2.14) holds and an easy computation shows
that $\mathbb{E}\|X\| \ln (1 + \|X\|) < \infty$. So by Lemma 5.5, Lemma 5.6,
Theorem 2.3, and Theorem 2.2 we obtain the following result.

\vskip 0.3cm

\begin{corollary}
Let $\mathbf{B}$ be a Banach space of stable type $1$. Let
$\{X_{n};~n \geq 1 \}$ be a sequence of independent copies of a
$\mathbf{B}$-valued random variable $X$. Then we have
\[
\sum_{n=1}^{\infty}\frac{1}{n}\left(\frac{\mathbb{E}\left\|S_{n}\right\|}{n}\right)
< \infty~~\mbox{if and only if}~~\mathbb{E}X = 0
~~\mbox{and}~~\int_{0}^{\infty}G(t) \ln \frac{1}{G(t)} dt < \infty.
\]
\end{corollary}

\vskip 0.3cm

For illustrating the conditions (2.11), (2.12), and (2.13) of
Theorem 2.3, we now present the following three examples.

\vskip 0.3cm

\begin{example}
Let $X$ be a real-valued random variable such that
\[
\mathbb{P}\left(X = - \frac{1}{1-a} \right) = 1 -a~~ \mbox{and}~~
\mathbb{P}(X > x) = \int_{x}^{\infty} \frac{1}{t^{2}\ln^{2}t} dt,
~~x \geq e
\]
where $a = \int_{e}^{\infty} \frac{1}{t^{2}\ln^{2}t} dt$. Then
\[
\mathbb{E}X = 0, ~~\mathbb{E}|X|\ln^{\delta}(1 + |X|) <
\infty~~\mbox{for all}~0 < \delta < 1,
\]
and, for all sufficiently large $n$,
\[
\mathbb{E}XI\{|X| \leq  n\} = - \mathbb{E}XI\{|X| > n\} = -
\int_{n}^{\infty} \frac{1}{t \ln^{2}t}dt = -\frac{1}{\ln n}.
\]
Note that
\[
\sum_{n=2}^{\infty}\frac{1}{n\ln n} = \infty.
\]
Thus (2.12) fails but, by Lemma 5.6, (2.13) holds.
\end{example}

\vskip 0.3cm

\begin{example}
Let $X$ be a real-valued symmetric random variable with density
function
\[
f(x) = \frac{b}{x^{2}(\ln|x|) (\ln\ln|x|)^{2}}I\{|x| > 3\},
\]
where $0 < b < \infty$ is such that $\int_{-\infty}^{\infty}f(x)dx =
1$. Clearly, (2.11) and (2.12) hold. Since
\[
\mathbb{P}(|X| > x) \sim \frac{2b}{x (\ln x) (\ln \ln
x)^{2}}~~\mbox{as}~ x \rightarrow \infty,
\]
we see that
\[
u_{n} \sim \frac{2bn}{(\ln n) (\ln \ln n)^{2}} ~~\mbox{as}~ n
\rightarrow \infty
\]
and hence, for all sufficiently large $n$,
\[
\begin{array}{lll}
\mbox{$\displaystyle \int_{u_{n}}^{n} \mathbb{P}(|X| > t) dt$}
&\geq& \mbox{$\displaystyle \int_{\frac{bn}{(\ln n) (\ln \ln
n)^{2}}}^{n} \frac{b}{t (\ln t) (\ln \ln t)^{2}}dt$}\\
&&\\
&\geq& \mbox{$\displaystyle \frac{b}{(\ln n)(\ln \ln n)^{2}}
\int_{\frac{bn}{(\ln n) (\ln \ln
n)^{2}}}^{n} \frac{1}{t} dt$}\\
&&\\
&\sim& \mbox{$\displaystyle \frac{b}{(\ln n) (\ln \ln
n)}~~\mbox{as}~ n \rightarrow \infty.$}
\end{array}
\]
Note that
\[
\sum_{n=3}^{\infty}\frac{b}{n(\ln n)(\ln \ln n)} = \infty
\]
and so (2.13) fails.
\end{example}

\vskip 0.3cm

\begin{example}
Let $X$ be a real-valued symmetric random variable with density
function
\[
f(x) = \frac{1}{2x^{2}}I\{|x| > 1 \}.
\]
Clearly, (2.12) holds. Since
\[
\mathbb{P}(|X| > x) = \frac{1}{x}, ~~x > 1,
\]
we have that
\[
u_{n} = n, ~~n \geq 1.
\]
Thus (2.13) also holds. However, (2.11) fails.
\end{example}

\vskip 0.5cm

\noindent
{\bf Acknowledgments}\\

\noindent The authors are extremely grateful to the Referee for very
carefully reading the manuscript and for offering numerous comments
and suggestions which enabled us to substantially improve the paper.
The Referee did not only provide comments to help us to improve the
presentation but more significantly the Referee presented improved
versions of some of the main results complete with their
proofs. Specifically, the Referee significantly improved our original
version of Theorem 3.1 and this improved version led to more elegant
proofs of many other results in the paper as was pointed out to us
by the Referee. The authors also thank Professor Andrzej Korzeniowski
for his interest in our work and for pointing out to us the relevance
of his paper [9] to ours. The research of Deli Li was partially supported
by a grant from the Natural Sciences and Engineering Research Council of
Canada and the research of Yongcheng Qi was partially supported by
NSF Grant DMS-0604176.

\newpage

{\bf References}

\begin{enumerate}

\item Azlarov, T. A., Volodin, N. A.: Laws of large numbers for
identically distributed Banach-space valued random variables. Teor.
Veroyatnost. i Primenen. {\bf 26}, ~584-590 (1981), in Russian.
English translation in Theory Probab. Appl. {\bf 26}, 573-580
(1981).

\item
Chow, Y.S., Teicher, H.: Probability Theory: Independence,
Interchangeability, Martingales, 3rd ed. Springer-Verlag, New York
(1997).

\item de Acosta, A.: Inequalities for {\it B}-valued random vectors
with applications to the law of large numbers. Ann. Probab. {\bf 9},
157-161 (1981).

\item Hechner, F., Heinkel, B.: The Marcinkiewicz-Zygmund LLN in
Banach spaces: A generalized martingale approach. J. Theor. Probab.
{\bf 23}, 509-522 (2010).

\item Hoffmann-J{\o}rgensen, J.: Sums of independent Banach
space valued random variables. Studia Math. {\bf 52}, 159-186
(1974).

\item Hoffmann-J{\o}rgensen, J., Pisier, G.: The law of large
numbers and the central limit theorem in Banach spaces. Ann. Probab.
{\bf 4}, 587-599 (1976).

\item
Klass, M. J.: Properties of optimal extended-valued stopping rules
for $S_{n}/n$. Ann. Probab. {\bf 1}~~719-757 (1973).

\item Kolmogoroff, A.: Sur la loi forte des grands nombres. C. R.
Acad. Sci. Paris S\'{e}r. Math. {\bf 191}, 910-912 (1930).

\item Korzeniowski, A.: On Marcinkiewicz SLLN in Banach spaces. Ann. Probab.
{\bf 12}, 279-280 (1984).

\item Ledoux, M., Talagrand, M.:  Probability in Banach Spaces:
Isoperimetry and Processes. Springer-Verlag, Berlin (1991).

\item Marcinkiewicz, J., Zygmund, A.: Sur les fonctions
ind\'{e}pendantes. Fund. Math. {\bf 29}, 60-90 (1937).

\item Maurey, B., Pisier, G.: S\'{e}ries de variables al\'{e}atoires
vectorielles ind\'{e}pendantes et propri\'{e}t\'{e}s
g\'{e}om\'{e}triques des espaces de Banach. Studia Math. {\bf 58},
45-90 (1976).

\item Marcus, M. B., Woyczy\'{n}ski, W. A.: Stable measures and
central limit theorems in spaces of stable type. Trans. Amer. Math.
Soc. {\bf 251}, 71-102 (1979).

\item Mourier, E.: El\'{e}ments al\'{e}atoires dans un espace de
Banach. Ann. Inst. H. Poincar\'{e} {\bf 13}, 161-244 (1953).

\item Pisier, G.: Probabilistic methods in the geometry of Banach
spaces, in Probability and Analysis, Lectures given at the 1st 1985
Session of the Centro Internazionale Matematico Estivo (C.I.M.E.),
Lecture Notes in Mathematics, Vol. {\bf 1206}, 167-241,
Springer-Verlag, Berlin (1986).

\item Rosi\'{n}ski, J.: Remarks on Banach spaces of stable type.
Probab. Math. Statist. {\bf 1}, 67-71 (1980).

\item Taylor, R. L.: Stochastic Convergence of Weighted Sums of
Random Elements in Linear Spaces, Lecture Notes in Mathematics, Vol.
{\bf 672}, Springer-Verlag, Berlin (1978).

\item Woyczy\'{n}ski, W. A.: Geometry and martingales in Banach
spaces-Part II: Independent increments, in Probability on Banach
Spaces (Edited by J. Kuelbs), Advances in Probability and Related
Topics Vol. {\bf 4} (Edited by P. Ney), 267-517, Marcel Dekker, New
York (1978).

\end{enumerate}

\end{document}